\documentclass[a4paper]{article}
\usepackage{amsmath, amstext, amsxtra, amsthm, amscd, amsgen,amsbsy, amsopn, amsfonts, latexsym, amssymb, amscd, amsthm}
\usepackage{marvosym}
\usepackage{color}

\def\wrt{with respect to}

\def\proofend{\hbox to 1em{\hss}\hfill $\blacksquare $\bigskip }

\newtheorem*{theorem*}{Main Theorem}

\newtheorem{theorem}{Theorem}[section]

\newtheorem{remark}[theorem]{Remark}

\def\Z{{\mathbb Z}}

\def\Q{{\mathbb Q}}
\def\C{{\mathbb C}}
\def\N{{\mathbb N}}

\begin{document}


\title{On the moduli space of nonnegatively curved metrics on Milnor spheres}

\author{Anand Dessai\thanks{The author acknowledges support by SNF grant 200021E-172469 and the DFG priority programme SPP 2026.}}

\date{}

\maketitle
\begin{abstract}{Let $M$ be a Milnor sphere or, more generally, the total space of a linear $S^3$-bundle over $S^4$ with $H^4(M;\Q )=0$. We
show that the moduli space of metrics of nonnegative sectional curvature on $M$ has infinitely many path components. The same holds true for the moduli space of metrics of positive Ricci curvature on $M$.}\end{abstract}

\noindent
\section{Introduction}
For a closed smooth manifold $M$, let ${\cal R}(M)$ be the set of Riemannian metrics on $M$, equipped with the $C^\infty $-topology, and let ${\cal R}_{sec \ge 0}(M)$ denote the subspace of metrics of nonnegative sectional curvature. The diffeomorphism group ${\rm Diff}(M)$ acts on  ${\cal R}_{sec \ge 0}(M)$ by pulling back metrics. The {\em moduli space} ${\cal M}_{sec \ge 0}(M):={\cal R}_{sec \ge 0}(M)/\hbox{Diff}(M)$ is the quotient space, equipped with the quotient topology. Corresponding notations will be used for other curvature bounds.

In this note we consider total spaces of linear $S^3$-bundles over $S^4$. These manifolds are very interesting, both from a topological and a geometric point of view. Milnor showed that among them there are manifolds which are homeomorphic to $S^7$ but not diffeomorphic to $S^7$ with its standard smooth structure \cite{Mi56}. Following modern terminology a homotopy sphere which is the total space of a linear $S^3$-bundle over $S^4$ will be called {\em Milnor sphere}.

Grove and Ziller showed that every Milnor sphere admits a metric of nonnegative sectional curvature and asked whether its moduli space has infinitely many components (see \cite[Thm. A and Problem 5.6]{GZ00}). As pointed out in \cite[p.~15]{De17} this question has an affirmative answer. The purpose of this note is to provide a detailed proof.

\begin{theorem*}\label{main theorem} For every Milnor sphere $M$ the moduli space ${\cal M}_{sec \ge 0}(M)$ of nonnegatively curved metrics has infinitely many path components.
\end{theorem*}

The main ingredients of the proof are the construction of metrics of nonnegative sectional curvature on certain cohomogeneity one manifolds due to Grove and Ziller, the diffeomorphism classification of Milnor spheres and an idea of Gromov and Lawson for distinguishing components of metrics of positive scalar curvature on compact spin manifolds. The proof is in a sense rather simple as it does not use more involved techniques like the APS-index theorem or the Kreck-Stolz invariant (which could also be used to give slightly different proofs for the results of this paper). 
At the end of this note we explain how the proof can be adjusted to show the analogous statement for the total space $M$ of every linear $S^3$-bundle over $S^4$ with $H^4(M;\Q )=0$.

Wolfgang Ziller has pointed out to me the recent work of Goodman \cite[v1]{Go17} in which the author used the Kreck-Stolz invariant to show that the moduli space of metrics of nonnegative sectional curvature has infinitely many path components for a certain family of $S^1$-bundles over $S^2$-bundles over $\C P^2$. More recently, Goodman has extended his results to cover also the Milnor spheres and the above $S^3$-bundles over $S^4$, again using the Kreck-Stolz invariant \cite[v2]{Go17}.

Prior to the work of Grove and Ziller only one exotic $7$-dimensional sphere was known to admit a metric of nonnegative sectional curvature, the so-called Gromoll-Meyer sphere \cite{GM74}. Among the $14$ exotic $7$-dimensional spheres $10$ can be described as Milnor spheres (see \cite{Mi59,EK62}) and, hence, each of these carries metrics of nonnegative sectional curvature by \cite{GZ00}. Recently, Goette, Kerin and Shankar \cite{GKS17} have put out a preprint in which the authors give a construction of metrics of nonnegative sectional curvature on exotic $7$-dimensional spheres including the remaining four cases. We don't know whether the Main Theorem extends to the latter ones.

The proof of the Main Theorem also shows that for every Milnor sphere $M$ the moduli space ${\cal M}_{Ric > 0}(M)$ of metrics of positive Ricci curvature has infinitely many connected components. This was known before. In fact, Wraith \cite{Wr11} has shown by different methods that for every homotopy sphere of dimension $4k-1$, $k\geq 2$, which bounds a parallelizable manifold the moduli space of metrics of positive Ricci curvature has infinitely many connected components (the corresponding statement for positive scalar curvature is classical, see \cite[pp. 131--132]{GL83}).

The study of the space of metrics under curvature bounds and its moduli space (or its observer moduli space) has been an active research area in the last years. In this note our focus is on Milnor spheres and on total spaces of linear $S^3$-bundles over $S^4$. For a survey of some of the many other recent developments and references we refer to \cite{WT15}.

This note is structured as follows. In the next section we review the Grove-Ziller construction \cite{GZ00} of nonnegatively curved metrics on principal $Spin(4)$-bundles, principal $SO(4)$-bundles and associated sphere bundles over $S^4$. This leads to infinitely many metrics on every Milnor sphere $M$ which are simultaneously of nonnegative sectional and positive scalar curvature and which can be extended to nonnegatively curved metrics on the associated disk bundles. The section also contains a review of the classification of Milnor spheres up to orientation preserving diffeomorphism via the $\mu$-invariant of Eells and Kuiper.

In the following section we use an idea of Gromov and Lawson \cite{GL83} to show that these metrics belong to pairwise distinct path components of the moduli space ${\cal M}_{sec \ge 0}(M)$. The argument also shows that these metrics evolve under the Ricci flow to metrics which belong to pairwise distinct components of ${\cal M}_{Ric> 0}(M)$.

In the last section we discuss the modifications needed to extend the result to the total space $M$ of a linear $S^3$-bundle over $S^4$ with $H^4(M;\Q )=0$.

\bigskip

{\em Acknowledgement.} Some of the results of this note were obtained and announced during the Mini-Workshop {\em Spaces and Moduli Spaces of Riemannian
Metrics} at the Mathematical Research Institute of Oberwolfach in January 2017  (see \cite{De17}). It is my pleasure to thank the institute, the organizers and the participants for the stimulating working environment. I also like to thank David Gonz\'alez-\'Alvaro and Michael Wiemeler for helpful comments on an earlier version of this note.

\section{Nonnegative curvature on Milnor spheres} Let $P$ be a closed connected smooth manifold with finite fundamental group on which a compact Lie group $G$ acts smoothly with cohomogeneity one. Then the quotient space is an interval and $P$ can be described as the union of disk bundles over the two nonprincipal orbits glued together along their boundaries.

In the case where the nonprincipal orbits are of codimension two Grove and Ziller constructed an invariant metric of nonnegative sectional curvature on each disk bundle such that the metric is a product metric near the boundary and the restriction to the boundary is isometric to the normal homogeneous space $G/H$ (\wrt \ a fixed bi-invariant metric on $G$). Here $H$ denotes the principal isotropy subgroup of the $G$-action. By gluing the two disk bundles isometrically and equivariantly along the boundary the cohomogeneity one manifold $P$ inherits a $G$-invariant metric of nonnegative sectional curvature  (see \cite{GZ00} for details on the construction of these metrics and for basic facts about cohomogeneity one manifolds). The construction involves various choices. In particular, it involves a rotationally symmetric nonnegatively curved metric on the $2$-disk $D^2$ which is of product form on a neighborhood of the boundary  (see \cite[Remark 2.7]{GZ00}) and which can be chosen to be of positive curvature on its complement.

In the same paper Grove and Ziller showed (among many other things) that the total space $P$ of any principal $(S^3\times S^3)$-bundle $P\to S^4$ admits a cohomogeneity one action by $S^3\times S^3\times S^3$ with finite principal isotropy subgroup, one dimensional nonprincipal isotropy subgroups and singular orbits of codimension two.

Here $S^3\times S^3\times 1$ acts freely on the fibers  via the principal action and $1\times 1 \times S^3$ covers the well-known cohomogeneity one action of $SO(3)$ on $S^4$ induced from the action of $SO(3)$ by conjugation on the space of symmetric real $3\times 3$-matrices with vanishing trace (see \cite[\S 3]{GZ00} for details).

We fix an invariant Grove-Ziller metric of nonnegative sectional curvature on $P$ as above. Note that $Spin(4)\cong S^3\times S^3\times 1$ acts freely and isometrically on $P$ and that the quotient $P^*:=P/(-1,-1)$ is the total space of the associated principal $SO(4)$-bundle over $S^4$. 

Let $M:=(P\times S^3)/(S^3\times S^3\times 1)$ denote the 7-dimensional manifold obtained from the Riemannian product of $P$ and the round sphere $S^3$ of radius one by dividing out the diagonal action of $S^3\times S^3\times 1$ (where the latter acts via the standard 4-dimensional $SO(4)$-representation on $S^3$).

Since the Riemannian product $P\times S^3$ has nonnegative sectional curvature it follows from the Gray-O'Neill formula \cite{ON66, Gr67} that the submersion metric on $M$ has nonnegative sectional curvature as well. Moreover, it is not difficult to see that every horizontal space $\cal H$ of the Riemannian submersion $P\times S^3\to M$ contains a plane of positive sectional curvature. In fact, one can always find a plane ${\cal E}\subset H$ which maps onto a tangent plane of $S^3$ under the projection $P\times S^3\to S^3$. Using the Riemannian product structure one computes directly that the sectional curvature of $\cal E$ is positive. Hence, again by the Gray-O'Neill formula the submersion metric on $M$ has positive scalar curvature everywhere (one can show that the same is already true for the metric on $P$).  

The Milnor sphere $M$ can also be described in terms of the principal $SO(4)$-bundle $P^*\to S^4$ as the total space of the sphere bundle of the associated oriented $4$-dimensional vector bundle $E\to S^4$. Sphere bundles of this kind will be called {\em linear} $S^3$-bundles over $S^4$.

The isomorphism classes of principal $SO(4)$-bundles over $S^4$ or, equivalently, of oriented $4$-dimensio\-nal vector bundles $E\to S^4$ are classified by $\pi _3(SO(4))\cong \Z \times \Z $. Up to isomorphism $E\to S^4$ is determined by its Euler class $e(E)$ and its first Pontrjagin class $p_1(E)$ and any pair $(e,p_1)$ with $p_1\equiv 2\cdot e\bmod 4$ can occur. If $e(E)$ is a generator of $H^4(S^4;\Z)$ then the total space $S(E)$ of the sphere bundle of $E$ is a homotopy sphere and is, in fact, homeomorphic to $S^7$ as first shown by Milnor using Morse theory \cite{Mi56}.

We fix a generator $u$ of $H^4(S^4;\Z)$ and denote by $P^*_k\to S^4$ the principal $SO(4)$-bundle with Euler class equal to $u$ and first Pontrjagin class equal to $2k\cdot u$, $k$ an odd integer. Let $E_k\to S^4$ be the associated vector bundle and let $M_k:=S(E_k)$ and $W_k:=D(E_k)$ be the total space of the associated unit sphere and disk bundle, respectively.

Recall that the Milnor sphere $M_k=P^*_{k} \times _{SO(4)} S^3$ comes with a Grove-Ziller metric of nonnegative sectional and positive scalar curvature. We equip the closed $4$-dimensional disk $D^4$ with a rotationally symmetric metric of nonnegative sectional curvature which is of product form near the boundary and which extends the round metric on the unit sphere. Now $W_{k}$ can be identified with the quotient $P^*_{k} \times _{SO(4)} D^4$ and, hence, inherits a submersion metric of nonnegative sectional curvature by the Gray-O'Neill formula \cite{ON66, Gr67}. For further reference we summarize the discussion in the following (see \cite{GZ00})

\begin{theorem}\label{nnc metric}
Every Milnor sphere $M_{k}$ admits a metric which is simultaneously of nonnegative sectional curvature and  positive scalar curvature and which extends to a metric of nonnegative sectional curvature on the disk bundle $W_{k}$ which is of product form near the boundary.\proofend
\end{theorem}

Next we review the classification of Milnor spheres up to orientation preserving diffeomorphism via the $\mu$-invariant of Eells and Kuiper.  We fix an orientation of $W_k$ by requiring that the square of a generator of $H^4(W_k,M_k;\Z )\cong \Z $ evaluated on the fundamental class is positive and orient $M_k$ accordingly. Note that $W_k$ is not parallelizable since its signature is equal to one. Note also that $M_{k}$ and $W_{k}$ have unique spin structures and that $M_{k}$ is the spin boundary of $W_{k}$.

Milnor showed that among the $7$-dimensional manifolds $M_k$, $k$ odd, some are not diffeomorphic to the standard sphere, thereby exhibiting the first examples of exotic 7-dimensional spheres \cite{Mi56}. To distinguish from the standard sphere Milnor first used an invariant $\lambda \in \Z/7\Z$ which is based on Hirzebruch's signature theorem and later introduced an invariant $\lambda ^\prime \in \Z/28\Z$ defined in terms of a parallelizable manifold with boundary $M_k$ \cite{Mi59}. Subsequently, Eells and Kuiper modified this approach and defined an invariant $\mu\in \Z/28\Z$ using spin manifolds, the signature and the $\hat A$-genus. All these invariants depend only on the oriented diffeomorphism type of $M_k$, are additive with respect to connected sum operation and change sign if one changes the orientation.

By the work of Milnor \cite{Mi59} and the $h$-cobordism theorem of Smale \cite{Sm62} the set of diffeomorphism classes of oriented closed smooth manifolds homeomorphic to $S^7$ form a cyclic group of order $28$ \wrt \ connected sum operation. Moreover, the elements of this group are determined by their $\lambda ^\prime $-invariant as well as by their $\mu$-invariant. In the following we will consider the Eells-Kuiper invariant computed in terms of $W_k$.
 
Let $W_k$ and $M_k$ be oriented as before. For $k=2h-1$ odd the invariant $\mu (M_k)\in  \Z/28\Z\subset \Q /\Z$ can be computed from $W_k$ to be  (see \cite[Eq. (11)]{EK62})
$$\mu (M_k)\equiv (p_1^2[W_k]-4\cdot \mathrm{sign}(W_k))/(2^7\cdot 7)\equiv \tilde h/28  \bmod \Z ,$$
where $\tilde h:=h(h-1)/2\in \Z $. Here $p_1^2[W_k]$ is defined by first taking a preimage $x$ of $p_1(W_k)$ under $H^4(W_k,M_k;\Q)\to H^4(W_k;\Q )$ and then evaluating $x^2$ on the fundamental class of $W_k$. The set of $\mu$-invariants of $M_k$, $k$ odd, can be computed to be equal to
$$\{r/28\in \Q /\Z\,\mid \, r=0,1,3,6,7,8,10,13,14,15,17,20,21,22,24,27\}.$$
Note that any of these elements in $\Q /\Z$ can be obtained from infinitely many distinct odd integers $k$. For further reference  we summarize the discussion in the following (see \cite{EK62})

\begin{theorem}\label{diffeo typ} The $28$ distinct closed oriented smooth manifolds homeomorphic to $S^7$ can be distinguished by their $\mu$-invariant. Among them $16$ are diffeomorphic to a Milnor sphere and any Milnor sphere is orientation preserving diffeomorphic to $M_k$ for infinitely many odd integers $k$.\proofend
\end{theorem}

 If $\mu (M_k)\equiv r/28 \bmod \Z $ with $14<r<27$ then $\mu(-M_k)$ is congruent to $ r/28 \bmod \Z $ for some $0< r < 14$. Therefore the $\mu$-invariants of the Milnor spheres (up to orientation) belong to
$$\{\pm r/28\in \Q /\Z\,\mid \, r=0,1,3,4,6,7,8,10,11,13,14\}.$$
Hence, ignoring orientation, $11$ of the $15$ distinct closed smooth manifolds homeomorphic to $S^7$  are diffeomorphic to a Milnor sphere.

We fix a Milnor sphere $M_l$ for some $l\geq 0$ and choose an infinite family of Milnor spheres $\{M_{k}\}_k$, such that each $M_k$ is orientation preserving diffeomorphic to $M_l$ and such that all $k$ are $\geq 0$. For example one can choose the family $M_{k}$, $k \in l +112\cdot \N $ (since $\mu(M_{l})=\mu(M_{l+112})$ in $\Q /\Z $ the Milnor spheres $M_{l}$ and $M_k$ are orientation preserving diffeomorphic for all $k \in l +112\cdot \N $). 

Let $g_k$ be the Grove-Ziller metric on $M_k$ considered before (see Thm. \ref{nnc metric}) and let $[g_k]$ denote the corresponding element in the moduli space ${\cal M}_{sec \ge 0}(M_l)$ obtained by pulling back $g_k$ via a diffeomorphism $M_l\to M_k$. Note that $[g_k]$ does not depend on the choice of the diffeomorphism and that we may choose the diffeomorphism to be orientation preserving. Also note that $[g_k]$ can be represented by a metric on $M_l$ which is simultaneously of nonnegative sectional and positive scalar curvature.

In the next section we will show that the elements $[g_k]$, $k\in  l +112\cdot \N$, belong to different path components of ${\cal M}_{sec \ge 0}(M_l)$.

\section{Proof of the Main Theorem}
We fix a Milnor sphere $M:=M_l$ and consider the elements $[g_k]\in {\cal M}_{sec \ge 0}(M)$, $k\in l+ 112 \cdot \N$, from the last section. We will show below that these elements belong to pairwise distinct path components of ${\cal M}_{sec \ge 0}(M)$. The proof is by contradiction, the idea may be summarized as follows: Consider $[g_{k_0}]$ and $[g_{k_1}]$,  $k_0\neq k_1$.

In the first step we show that a path in ${\cal M}_{sec \ge 0}(M)$ connecting $[g_{k_0}]$ and $[g_{k_1}]$  can be lifted and deformed to a path $\gamma:[0,1]\to {\cal R}_{scal> 0}(M)$ in the space of metrics of positive scalar curvature such that $\gamma (0)$ represents $[g_{k_0}]$, $\gamma (1)$ represents $[g_{k_1}]$ and such that $\gamma (t)$ is constant near the end points. An elegant way to obtain the deformation from nonnegative sectional curvature to positive scalar curvature (and even positive Ricci curvature) uses work of B\"ohm and Wilking on the Ricci flow \cite{BW07}.

After stretching the interval we may assume that the metric on $M\times [0,a]$, $a\gg 0$, defined by the metrics $\widetilde \gamma (t):=\gamma (t/a)$, $t\in [0,a]$, has positive scalar curvature.

In the second step we use a variant of the following elementary version of the relative index of Gromov and Lawson for compact spin manifolds (see \cite{GL83}):

Suppose $g_0$ and $g_1$ are two metrics of positive scalar curvature on a closed spin manifold $N$ of dimension $4k-1$. Suppose in addition that $g_i$ extends to a metric $h_i$ of nonnegative scalar curvature on a spin manifold $Y_i$ with spin boundary $N$ and that $h_i$ is of product form near the boundary. Then the $\hat A$-genus of the closed $4k$-dimensional spin manifold $Y_0\cup _N (-Y_1)$ vanishes if $g_0$ and $g_1$ belong to the same component of ${\cal R}_{scal> 0}(N)$.

We apply this to our Milnor manifold and to the disk bundles considered in Theorem \ref{nnc metric}.

In the last step we compute the $\hat A$-genus of the corresponding closed $4k$-dimensional spin manifold to be non-zero, thereby arriving at the desired contradiction. Here are the details. 

\begin{theorem}[Main Theorem] For $k_0, k_1\in l+ 112 \cdot \N$,  $k_0\neq k_1$, the classes $[g_{k_0}]$ and $[g_{k_1}]$ belong to different path components of ${\cal M}_{sec \ge 0}(M)$.
\end{theorem}

\bigskip
\noindent
{\bf Proof:} Let us assume to the contrary that $[g_{k_0}]$ and $[g_{k_1}]$ belong to the same path component. We fix orientation preserving diffeomorphisms $\Psi_i:M\to M_{k_i}$ for $i=0,1$. Recall from the end of the last section that the element $[g_{k_i}]\in {\cal M}_{sec \ge 0}(M)$ is represented by the metric $\Psi_i^*(g_{k_i})$.

Now consider a path connecting $[g_{k_0}]$ and $[g_{k_1}]$ in $ {\cal M}_{sec \ge 0}(M)$. By Ebin's slice theorem \cite{Eb70} the path can be lifted to a path $\hat \gamma (t)$, $t\in I:=[0,1]$, in $ {\cal R}_{sec \ge 0}(M)$ which starts in $\hat \gamma (0)=\Psi_0^*(g_{k_0})$. Hence, the Riemannian manifolds $(M,\hat \gamma (0))$ and $(M_{k_0},g_{k_0})$ are isometric by the orientation preserving diffeomorphism $\Psi_0$. The metric $\hat \gamma (1)$ lies in the same ${\rm Diff}(M)$-orbit of $\Psi_1^*(g_{k_1})$. If $\hat \gamma (1)$ is obtained as the pullback of $\Psi_1^*(g_{k_1})$ by an orientation {\em reversing} diffeomorphism of $M$ we will replace $g_{k_1}$ by its pullback under an orientation {\em reversing} diffeomorphism of $M_{k_1}$. So in the following we can assume that $(M,\hat \gamma (i))$ and $(M_{k_i},g_{k_i})$ are isometric by an orientation preserving diffeomorphism for $i=0,1$.

As shown by B\"ohm and Wilking \cite{BW07} the metric $\hat \gamma (t)$ for fixed $t$ evolves under the Ricci flow instantly to a metric on $M$ of positive Ricci curvature. Hence, the path $\hat \gamma (t)$, $t\in I$, evolves under this flow immediately to a path in the space of metrics of positive Ricci curvature. Concatenation of the evolved path and the trajectories of the end points of $\hat \gamma (t)$ then yields a path $\gamma $ connecting the end points $\hat \gamma (0)$ and $\hat \gamma (1)$. Note that $\gamma$ is a path in the space of metrics of positive scalar curvature, ${\cal R}_{scal> 0}(M)$, since the metric $\gamma (t)$ has positive scalar curvature at the end points of the interval and positive Ricci curvature for points in the interior. After reparametrization and a small perturbation leaving the endpoints fixed we may assume that $\gamma :I\to {\cal R}_{scal> 0}(M)$ is a smooth path which is constant close to the end points $\gamma (0)=\hat \gamma (0)$ and $\gamma (1)=\hat \gamma (1)$.

The metrics $\gamma (t)$, $t\in I$, define a metric on $M\times I$ which is fiberwise of positive scalar curvature. After stretching the interval $I$ sufficiently (i.e. after increasing the distance between different fibers sufficiently) one obtains a metric on $M\times [0,a]$, $a\gg 0$, which has positive scalar curvature (see \cite[Lemma 3]{GL80}). Let $\widetilde \gamma (t):=\gamma (t/a)$ denote the corresponding path in ${\cal R}_{scal> 0}(M)$.

Next we will use the elementary relative index invariant of Gromov and Lawson for compact spin manfolds (see \cite[pp. 130--131]{GL83}) to derive a contradiction. We fix orientation preserving isometries $\varphi _i:(M_{k_i},g_{k_i}) \to (M,\hat \gamma (i))=(M,\widetilde \gamma (a\cdot i))$ for $i=0,1$. Recall from Theorem \ref{nnc metric} that $g _{k_i}$ extends to a metric, say $h_{k_i}$, of nonnegative sectional curvature on $W_{k_i}$ which is of product form near the boundary. By gluing $(W_{k_0}, h_{k_0})$, $-(W_{k_1}, h_{k_1})$ and $(M\times [0,a], \widetilde \gamma (t))$ via $\varphi _0:(M_{k_0},g_{k_0}) \to (M,\widetilde \gamma (0))$ and $\varphi _1:(M_{k_1},g_{k_1}) \to (M,\widetilde \gamma (a))$ together we obtain a closed $8$-dimensional Riemannian  spin manifold $$X=W_{k_0}\cup_{\varphi _0}(M\times [0,a])\cup_{\varphi _1}(-W_{k_1})$$ which is of nonnegative scalar curvature and has positive scalar curvature somewhere. By Lichnerowicz' theorem \cite{Li63} the Dirac operator of $X$ is invertible and its index $\hat A(X)$ vanishes.

Note that $X$ is diffeomorphic to $W_{k_0}$ and $-W_{k_1}$ glued together via an orientation preserving diffeomorphism. Now the signature of $X$,  $\mathrm{sign}(X)$, does not depend on the choice of the gluing map and is equal to the sum of the signatures of $W_{k_0}$ and $-W_{k_1}$ (see for example \cite{Mi56, Mi59}, \cite[p. 98]{EK62}). Hence,
$$\mathrm{sign}(X)=\mathrm{sign}(W_{k_0})-\mathrm{sign}(W_{k_1})=1-1=0.$$

The vanishing of the $\hat A$-genus and the signature of $X$ forces the Pontrjagin number $\langle p_1(X)^2,[X]\rangle$ to vanish as well (here $[X]$ denotes the fundamental class of $X$). The Pontrjagin number can be computed in terms of the disk bundles $W_{k_i}$. More precisely, one has
$$\langle p_1(X)^2,[X]\rangle = p_1^2[W_{k_0}] -p_1^2[W_{k_1}] ,$$
where, as before, $p_1^2[W_{k_i}]$ is defined by first taking a preimage of $p_1(W_{k_i})$ under $H^4(W_{k_i},M_{k_i};\Z)\to H^4(W_{k_i};\Z )$ and then evaluating its square on the fundamental class of $W_{k_i}$, $i=0,1$ (see \cite{EK62}).

The tangent bundle of the disk bundle $W_{k_i}=D(E_{k_i})$ is stably isomorphic to $\pi ^*(E_{k_i})$, where $\pi$ denotes the projection $W_{k_i}\to S^4$.

Hence, $p_1(W_{k_i})=2\cdot k_i\cdot \pi^*(u)\in H^4(W_{k_i};\Z )$, where $u\in H^4(S^4;\Z )$ is the fixed generator. This gives a contradiction since $p_1^2[W_{k_0}] -p_1^2[W_{k_1}]=(2k_0)^2-(2k_1)^2$ is non-zero. \proofend

\begin{remark}
The proof shows that the metrics obtained by applying the Ricci flow to lifts of $[g_k]$, $k\in l+112\cdot \N$, represent pairwise distinct components of the moduli space of metrics of positive Ricci curvature ${\cal M}_{Ric> 0}(M)$ (and ${\cal M}_{scal> 0}(M)$ as well). Hence, for any Milnor sphere $M$ the moduli space\linebreak  ${\cal M}_{Ric> 0}(M)$ has infinitely many connected components (cp. \cite{Wr11}).
\end{remark}

\begin{remark} In the recent preprint \cite{GKS17} Goette, Kerin and Shankar give a construction of metrics of nonnegative sectional curvature on exotic $7$-dimension\-al spheres including the non-Milnor spheres. The latter can be described as $S^3$-orbi-bundles but do not come with any obvious null-cobordism with nice geometric properties. This makes it difficult to treat their moduli spaces with the methods above.
\end{remark}

\section{Moduli spaces for linear $S^3$-bundles over $S^4$}

In this section we discuss the modifications needed to extend the Main Theorem to the total space $M$ of a linear $S^3$-bundle over $S^4$ with $H^4(M;\Q )=0$. The latter condition is fulfilled precisely when the Euler class of the underlying $4$-dimensional vector bundle is non-trivial. In the case where the Euler class is trivial the first Pontrjagin class of $M$ takes values in $H^4(M;\Z )\cong \Z $ and determines the vector bundle up to isomorphism. So the argument for the Milnor spheres does not apply to this situation. However, in the case where the Euler class is non-trivial the proof of the Main Theorem does extend with minor modifications.

\begin{theorem} Let $M$ be the total space of a linear $S^3$-bundle over $S^4$ with $H^4(M;\Q )=0$. Then the moduli space of metrics of nonnegative sectional curvature on $M$ has infinitely many path components. The same holds true for the moduli space of metrics of positive Ricci curvature on $M$.
\end{theorem}

\bigskip
\noindent
{\bf Proof:} Recall that the isomorphism classes of principal $SO(4)$-bundles over $S^4$, or, equivalently, the isomorphism classes of oriented $4$-dimensional real vector bundles $E\to S^4$ are determined by the Euler class $e(E)$ and the first Pontrjagin class $p_1(E)$ and any pair $(e,p_1)$ with $p_1\equiv 2\cdot e\bmod 4$ can occur. We restrict to bundles with non-trivial Euler class since $H^4(M;\Q )=0$. Let $u\in H^4(S^4;\Z )$ be a fixed generator. Then the Euler class of $E$ can be written as $n\cdot u$ for some integer $n\neq 0$. Note that $n$ changes sign if we change the orientation of the bundle. Since the moduli space of the associated sphere bundle does not depend on the orientation we may assume from now on that $n>0$.

For $k\in \Z$, $k\equiv n\bmod 2$, let $E_k$ denote the oriented $4$-dimensional real vector bundle over $S^4$ with $p_1(E_k)=2k\cdot u$ (and $e(E_k)=n\cdot u$). Let $M_k$ and $W_k$ denote the associated sphere and disk bundles, respectively. The integral cohomology of $M_k$ is concentrated in degree $0,4,7$, with  $H^0(M_k;\Z )\cong H^7(M_k;\Z )\cong \Z $ and $H^4(M_k;\Z)\cong \Z /n\Z $. As before we orient $W_k$ and $M_k$ such that $sign(W_k)=+1$.

Now the manifold $M$ in the theorem is diffeomorphic to $M_l$ for some $n>0$ and some $l\in \Z$ with $l\equiv n\bmod 2$. So our aim is to show that ${\cal M}_{sec\geq 0}(M_l)$ and ${\cal M}_{Ric> 0}(M_l)$  have infinitely many path components.

As in the proof of the Main Theorem we want to find for $M=M_l$ infinitely many $M_{k_i}$ of nonnegative sectional curvature such that each $M_{k_i}$ is diffeomorphic to $M_l$ as an oriented manifold and such that the induced metrics on $M_l$ represent pairwise distinct path components of the moduli space.

The classification of total spaces of linear $S^3$-bundles over $S^4$ up to homeomorphism was considered by Tamura already shortly after Milnor exhibited the first examples of exotic $7$-dimensional spheres. In particular, it is shown in \cite[Thm. 3.1]{Ta58}, that $M_k$ and $M_{k^\prime}$ are homeomorphic if $k^\prime\equiv k \bmod 2n$.

From smoothing theory one knows that the obstructions for being $PL$-homeomorphic or diffeomorphic belong to finite sets since $TOP/PL$ is an Eilen\-berg-MacLane space $K(\Z /2\Z,3)$ (see \cite{KS77}) and $PL/{O}$ is $6$-connected with $\pi _7(PL/{O})$ cyclic of order $28$ (see \cite{Mi59,Sm62,KM63,Ce68,Wa70,Br72,HM74}). This implies that the family $\{M_{k^\prime}\, \mid \, k^\prime\equiv k \bmod 2n \}$ contains an infinite subfamily of pairwise diffeomorphic manifolds. Moreover, it follows from the work of Sullivan \cite{Su77} that the entire family $\{M_{k}\}_k$ contains only finitely many diffeomorphism types.

In order to obtain more precise information we will use the classification of total spaces of linear $S^3$-bundles over $S^4$ with vanishing rational cohomology in degree $4$. These spaces (and other two-connected $7$-manifolds) have been classified to a large extent by Wilkens who applied handlebody theory building on techniques of Wall. Crowley and Escher completed the classification using the topological Eells-Kuiper invariant of Kreck and Stolz (see \cite{CE03} and references therein).

From this classification one knows that $M_k$ and $M_{k^\prime}$ are diffeomorphic as oriented manifolds if $k^\prime\equiv k \bmod 112n$ (see \cite[Cor. 1.6]{CE03}, where $M_{m,n}$ cor\-responds to $M_{n+2m}$ in our notation).

Now the proof follows the line of argument for the Main Theorem: For $M:=M_l$, $l\geq 0$ fixed, one considers the family $M_{k}$, $k \in l +112n\cdot \N $. Each $M_k$ is diffeomorphic to $M_l$ by an orientation preserving diffeomorphism and comes with a Grove-Ziller metric $g_k$ of nonnegative sectional curvature with the properties described in Theorem \ref{nnc metric}. By computing the relative index invariant one finds that the metrics $g_k$ induce elements in the moduli space ${\cal M}_{sec \ge 0}(M)$ which belong to pairwise distinct path components. Moreover, the corresponding metrics on $M$ evolve under the Ricci flow to metrics which represent pairwise distinct components of ${\cal M}_{Ric> 0}(M)$.  \proofend

\end{document}